\widowpenalty=15000
\overfullrule=0pt
\pretolerance=-1
\tolerance=2500
\doublehyphendemerits 50000
\finalhyphendemerits 25000
\adjdemerits 50000
\hbadness 1500
\abovedisplayskip=7pt plus 5pt minus 6pt
\belowdisplayskip=7pt plus 5pt minus 6pt
\hsize39pc \vsize52pc
\baselineskip11.5pt
\magnification=1200
\parindent 0pt
\def\exp{{\rm \hskip1.5pt  exp \hskip1.5pt }}
\def\summe{\sum\nolimits}
\font\gross=cmbx12
\font\grossrm=cmr12
\font\gr=cmr17

%\font\nineit=amti9
%\font\eightit=amti8

%
%\font\calligX=callig15 at 10pt
%\newfam\calligfam
%\textfont\calligfam=\calligX
%\font\calligVIII=callig15 at 8pt
%\scriptfont\calligfam=\calligVIII
%\font\calligV=callig15 at 5pt
%\scriptscriptfont\calligfam=\calligV
%\def\script{\fam\calligfam\calligX}
%
\def\parein#1#2{\par\noindent\rlap{\rm {#2}}\parindent#1\hang\indent\ignorespaces}
\def\paraus{\parindent 0pt\par}

\def\l{\ell}
\def \ol{\overline}
\def \1{\backslash}
%\input NEWMAKROS.tex
%
%
% output parameter 
%\input NEWMAKROS.tex
%\output={
%\ifodd\pageno\hoffset=0pt
%\else
%\hoffset=-0.27truecm\fi\plainoutput}
%
%
%
\def\erfc{{\rm Erfc}\,}
\def\re{{\rm Re}\,}

\def\la{{\cal L}}
\def\fa{{\cal F}}

\def\chyp #1 { C_{{\rm hyp}, #1 }}

\font\authorfont=cmcsc10 at 12pt
\font\titlefont=ptmb at 13pt%cmbx12 at 14pt 
\font\affiliationfont=cmti12 at 11pt
 
\font\gross=cmbx12

\font\ninecaps=cmcsc9
\font\ninebf=cmbx9
\font\nineit=cmti9
\font\ninerm=cmr9

\abovedisplayskip=8pt plus 5pt minus 3pt
\belowdisplayskip=8pt plus 5pt minus 3pt
\baselineskip=15pt
%
%\centerline{  } \pageno=0 \vfill\eject
%
\centerline{\titlefont ON THE VALUATION OF ARITHMETIC-AVERAGE ASIAN}
\centerline{\titlefont OPTIONS: LAGUERRE SERIES AND THETA INTEGRALS}
\baselineskip=11.5pt
\vskip.3cm
\centerline{\authorfont Michael  Schr\"oder}
\vskip.2cm
\centerline{\affiliationfont Lehrstuhl Mathematik III}
%\vskip.05cm
\centerline{\affiliationfont Seminargeb\"aude A5,
Universit\"at Mannheim,  D--68131 Mannheim}
\vskip.5cm 
\centerline{ 
\vbox{
\hsize=10.85cm
\baselineskip=9truept
\ninerm
In a recent significant advance, using Laguerre series, 
the valuation of Asian options has been reduced in [{\ninebf D}] 
to computing the negative moments of Yor's accumulation 
processes for which 
%Great efforts have been made by Yor and his 
%coworkers to analyze the probabilistic structure and probabilistic
%consequences of the 
functional recursion rules are given. Stressing the role of Theta 
functions, this paper now solves these recursion rules and 
expresses these negative moments as linear combinations of certain 
Theta integrals. Using the Jacobi transformation formula, very 
rapidly and very stably convergent series for them are derived. 
In this way a computable series for Black--Scholes
price of the Asian option results which is numerically illustrated. 
Moreover, the Laguerre series approach of [{\ninebf D}] is made 
rigorous, and extensions and modifications are discussed. 
The key for this is the analysis of the integrability and 
growth properties of the Asia density in [{\ninebf Y}], basic 
problems which seem to be addressed here for the first time.}}
\vskip.5cm
{\bf 1.\quad Introduction:}\quad Asian options are path--dependent 
options on the arithmetic average of the price of their underlying 
security. While they are widely traded 
financial securities their valuation still poses intriguing problems, 
even in the Black--Scholes setting. The Laguerre series approach  
of [{\bf D}] so was a significant advance and reduced valuing Asian 
options to computing the negative moments of the averaging process. 
For the latter it gave functional recurrence rules. Their  
structure and significance have been thoroughly analyzed by Yor 
and his coworkers from a probabilistic point of view. 
\medskip   
Taking an analytic point of view, this paper explains how to express
these negative moments using Theta functions, and in this way derives
computable Laguerre series for the value of the Asian option different
from those of [{\bf SE}]. The idea is that using the recurrence relations
of [{\bf D}] the negative higher moments of the averaging processes 
should be linear combinations of certain Theta integrals. This  
is true in the basic case of \S  10, and depending on the relative 
constellation of risk neutral drift and squared volaltility only 
finitely many correction terms have to be added for the general case 
of \S 12. Using the Jacobi transformation formula for Theta series,  
\S 11 derives a series for these Theta integrals which is optimized 
with respect to speed and stability of convergence. We so obtain a 
Laguerre series for the value of the Asian option whose coefficients 
are series given by integration against Theta series.
\medskip
This occurence of Theta functions in valuing the Asian option is 
rather surprising. Indeed, higher moments of the averaging 
process have been thoroughly studied in [{\bf YE}, Part A]. In 
their simplest form they were found to be expressed by Bougerol's 
identity as the expectation of the respective powers of the 
hyperbolic sine evaluated on a Brownian motion. While a structural 
explanation of this link of Asian options with modular is still 
missing, it highlights a characteristic difficulty of the Laguerre 
series approach. Its workability crucially depends on the specifics 
of the option to be valued.
\medskip
For the Asian option, we think any valid expansion into orthogonal 
functions will eventually be based on Yor's Asia density of 
[{\bf Y}, \S 6]. However, it is measurability which is 
addressed in [{\bf Y}, \S 6] while for such purposes integrability 
is required. Thus we establish  in \S 6 integrability of a certain 
class of functions of exponential type with respect to the Asia 
density, and put this result to work in deriving two Laguerre series
for the value of the Asian option. 
\medskip
\nopagenumbers
\def\Verfasser{{\ninecaps Michael Schr\"oder}}
\def\Kapiteltitel{{\ninecaps Laguerre series and Theta integrals for 
Asian options}}
\headline={\vbox{\line{\ninebf \ifodd\pageno\hss \Kapiteltitel\hss
\folio
\else
\folio\hss\Verfasser\hss\fi}
%\hrule
}\hss}
First, we illustrate a natural principle 
for getting series expansions of option prices in terms of higher 
moments of the option's control variable. The idea is to Laguerre 
expand in the risk neutral expectation that gives this value the 
taking--the--non--negative--part function, and 
then try to put the expectation through the series so obtained. 
In \S 7 we combine our integrability results with characteristic
mean convergence convergence results to justify
this last operation, and so make this series available to valuing 
Asian options. 
\medskip
Second, there are the Laguerre series in the spirit 
of [{\bf D}]. Their idea is to construct the price of the Asian 
option as probability density in the strike price using the notion 
of ladder height densities, Laguerre expand this density, and try to 
re--interpret the Laguerre coefficients in terms of the reciprocal 
of the averaging processes. As explained in \S 8, Laguerre 
expandability of the Asian ladder height densities crucially depends 
on the integrability results of \S 6, and we have slightly modified 
the original series in negative moments of the averaging process 
of [{\bf D}] in addition. 
\medskip
Moreover, we give a second way of re--interpreting the Laguerre
coefficients  using expectations of gamma function type 
integrands evaluated at reciprocals of the average process.  
These expectations also come up naturally if convergence 
of the Laguerre series in negative moments is studied. Indeed, 
empirical evidence in \S 13 suggests that their convergence behaviour 
is controlled by a certain convergence parameter. Characterizing
its optimal values will be in terms of these gamma type 
expectations. As indicated in \S 14, the latter now depend on 
Bessel functions, so have a different character than the 
negative moments, and we hope to return to their study elsewhere.  
\goodbreak
\bigskip
{\bf 2.\quad Preliminaries on Laguerre series:}\quad
This section collects pertinent properties of Laguerre polynomials from
[{\bf L},  \S 4] and [{\bf S}] fixing any real number $\alpha>-1$.
%\medskip
For any non--negative 
integer $n$, the {\it $\alpha$--Laguerre polynomial\/} 
$L_n^\alpha(z)$
%$L_n^{\raise-2pt\hbox{$\scriptstyle \alpha$}}$ 
is then given by
$$
L_n^\alpha(z)=\sum\nolimits_{k=0}^n (-1)^k 
{\Gamma(n\!+\!\alpha\!+\!1)\over\Gamma(k\!+\!\alpha\!+\!1)}
{z^k\over k! (n\!-\!k)!}\, , 
$$ 
for any complex number $z$. The first few $\alpha$--Laguerre polynomials
so are
$L_0^\alpha (z)=1$, $L_1^\alpha (z)=1\!+\!\alpha\!-\!z$, and
$L_2^\alpha (z)
=(1/2)\cdot( (1\!+\!\alpha)(2\!+\!\alpha)-2(2\!+\!\alpha)z+z^2)$,
and for any positive intger $n$ they satisfy the recurrence relation
$$
L_{n+1}^\alpha(z)= {2n\!+\!1\!+\!\alpha\!-\!z\over n\!+\!1}L_n^\alpha(z)
-{n\!+\!\alpha\over n\!+\!1}L_{n-1}^\alpha(z)\, . 
$$
%Specializing to $\alpha=0$, the {\it classical\/} Laguerre polynomials 
%$L_n$ are obtained which are explicitly given by
%$$
%L_n(z)=L_n^0(z)=\sum\nolimits_{k=0}^n 
%{(-z)^k\over k!}
%{n\choose k}\, .
%$$
The $\alpha$--Laguerre polynomials are orthogonal on the positive 
real line with respect to the weight $w_\alpha(x)=x^\alpha\exp(-x)$ 
in that  
$\int_{\raise1.5pt\hbox{$\scriptstyle 0$}}^\infty 
w_\alpha(x)(L_n^\alpha(x))^2L_m^\alpha(x)\, dx
=\Gamma(n\!+\!\alpha\!+\!1)/n! 
%={\Gamma(n\!+\!\alpha\!+\!1)\over n!}\, , 
$
and 
$
\int_{\raise1.5pt\hbox{$\scriptstyle 0$}}^\infty w_\alpha(x)
L_n^\alpha(x)L_m^\alpha(x)\, dx=0$
for any non--negative integers $n\ne m $. To discuss when functions $F$ 
on the positive real line have a Fourier--type expansion with respect to 
the $\alpha$--Laguerre polynomials, consider the case where $F$ is obtained 
by integrating up a function. So suppose there is a function $f$ on the 
positive real line which is integrable on any finite subinterval of the 
positive real line such that 
$F(x)=F(0)+\int_{\raise1.5pt\hbox{$\scriptstyle 0$}}^x f(y)\, dy$, 
for any $x>0$. Moreover assume that the 
functions $\sqrt{w_\alpha\, }F$ and $\sqrt{w_{\alpha+1}\, }f$ are square 
integrable on $(0,\infty)$. 
Then $F$ is represented by the {\it $\alpha$--Laguerre series\/} 
$$
e\vbox to 9pt{}^{-{\scriptstyle x\over\scriptstyle 2}}
x^{-{\scriptstyle x\over\scriptstyle 2}} F(x)
=
\sum_{n=0}^\infty 
c_n  e\vbox to 9pt{}^{-{\scriptstyle x\over\scriptstyle 2}}
x^{-{\scriptstyle x\over\scriptstyle 2}}
L_n^\alpha(x)
$$
which converges absolutely and uniformly for $x$ in any subinterval
$[c,\infty)$ of the positive real line. Any $n$--th
generalized Fourier coefficient $c_n$ of this series is 
given by
$$
c_n={n!\over \Gamma(n\!+\!\alpha\!+\!1)}
\int_0^\infty e^{-x}x^\alpha F(x) L_n^\alpha(x)\, dx\, . 
$$
Alternatively, if $F$ is continuous and $\sqrt{w_\alpha\, }F$ is 
square integrable on $(0,\infty)$, then $F$ is also  represented 
by the above Laguerre series, see [{\bf L}, \S 4.23]. 
Convergence at zero of this series, however, needs a separate study, see
for instance [{\bf S}, p.366f]. 
\goodbreak\bigskip
{\bf 3.\quad Preliminaries on Theta functions:}\quad  Following  
[{\bf M}], the concept to be discussed is the classical 
Riemann {\it Theta function} $\vartheta$ given for any complex number 
$z$ and any positive real number $t$ by:
$$ \vartheta(z|\, t)
={1\over \sqrt{ \pi t\, }}\summe _{n\in {\bf Z}}\ 
e\vbox to 9pt{}^{
-\, (z+n)^2\cdot 
{\scriptstyle1\over\scriptstyle t}
}.$$
This series converges absolutely, and uniformly on compact sets. 
Thus $\vartheta$ can be seen as a holomorphic function on the product of 
the complex plane with the upper complex half--plane. In modular forms 
this is usually done such that the above Theta function would 
be considered as evaluated not at $t$ but at the point $i\pi t$ of the 
upper half--plane. 
Theta functions are a basic class of holomorphic modular forms. They 
have been studied since the eightteenth century and are interrelated
with a number of areas of central importance for mathematics, as number 
theory, algebraic geometry, classical analysis, and partial differential 
equations. 
\medskip
Two such interrelations are to be described. For the first 
recall that $t^{-1/2}\exp(-\pi x^2/t)$ is the fundamental solution of 
the Heat equation on the line with initial data at $t=0$ a delta function 
at $x=0$. Thus $\vartheta$ at real arguments can be seen as the 
superposition of infinitely many  such solutions with initial data being 
delta functions at the half--integers $x=n+1/2$. 
\medskip
From the point of view 
of modular forms, Theta functions can be characterized by a certain 
periodicity behaviour with respect to each of their two variables. Here, 
the behaviour with respect to the second variable is deeper and more 
subtle. It is expressed by a functional equation with respect to 
the second variable, the {\it Jacobi transformation formula\/}. 
Restricting to real arguments, this formula is the following 
remarkable identity:
$$ 1+2\summe_{n=1}^\infty (-1)^n e^{-(\pi n )^2\cdot t}\cos(2\pi nx)
%=\vartheta(x|\,{t/\pi})
={1\over \sqrt{\pi t\, }}\summe _{n\in {\bf Z}} \ e\vbox to 9pt{}^{
-\, (z-{\scriptstyle 1\over\scriptstyle 2 }+n)^2\cdot 
{\scriptstyle1\over\scriptstyle t}
},$$
for any real numbers $x$ and  $t>0$. 
%[{\bf FB}, Kapitel VI \S 4, p.344].
The left hand side of this identity, the Jacobi transform of  
$\vartheta$, is rapidly converging for $t$ large, whereas its right hand 
side is rapidly converging for $t$ small. The Jacobi Transformation Formula 
is proved by Fourier analysis, exhibiting $\vartheta$ as the Fourier 
expansion with respect to the second variable of its Jacobi 
transform.
\goodbreak\bigskip
{\bf 4.\quad Basic notions:}\quad First we  discuss the notions 
basic for the analysis of Asian options. We work  in the Black--Scholes 
framework using the risk--neutral approach to the valuation of 
contingent claims.  
In this set--up there are two securities. There is a riskless security,
a bond, whose price grows at the continuously
compounding positive interest rate $r$. There is also a risky
security, whose price process $S$ is modelled as follows. Consider
a complete probability space equipped with the standard
filtration of a standard Brownian motion  on the time set
$[0,\infty)$. On this filtered space, we have the risk neutral
measure $Q$, which is a probability measure equivalent to the
given one. And then we have a standard $Q$--Brownian motion $B$ 
such that $S$  is the strong solution of the following stochastic 
differential equation:
$$ dS_t=\varpi\cdot S_t\cdot dt +\sigma \cdot S_t\cdot dB_t, 
\qquad t\in[0,\infty).$$  
The positive constant $\sigma$ is the volatility of $S$. The 
specific form of the otherwise arbitrary constant $\varpi$ depends 
on the nature of the security modelled (eg. stock, currency, commodity 
etc.). For example it is the interest rate if $S$ is a 
non--dividend--paying stock.
\medskip
Fix any time $t_0$ and consider the accumulation process $J$  given for 
any time $t$ by:
$$ J(t)=\int_{t_0}^t S_u \, du\, .$$
The European--style {\it arithmetic--average Asian option\/} written 
at time $t_0$, with maturity $T$, and fixed--strike price $K$ is then
the contingent claim on the closed time interval from $t_0$ to $T$ 
paying
$ (J(T)/(T\!-\!t_0)-K)^+:= \max\{ 0,J(T)/(T\!-\!t_0)-K\}$ 
%$$ \left( {J(T)\over T-t_0}-K\right)\vbox to 10pt {}^+=\max\left\{
%0, {J(T)\over T-t_0}-K\right\}$$
at time $T$. Recall that points in time are taken to be non--negative 
real numbers.
The price $C_t$ of the Asian option at any time $t$ between $t_0$ 
and $T$ is given as the following  risk neutral expectation 
$$C_t=e^{-r(T-t)} E^Q\left[\left( 
{J(T)\over T-t_0}-K\right)^+\Big| \fa_t\right]$$ 
which is conditional
on the information $\fa_t$ available at time $t$.
However, following [{\bf GY}, \S 3.2], do not focus on 
this price. As described there in great detail, we instead 
normalize the valuation problem, consider the factorization:
$$ C_t={e^{-r(T-t)}\over T-t_0}\cdot {4S_t\over \sigma^2}\cdot
C^{(\nu)}(h,q)\, ,  $$
and so reduce to computing 
$$ C^{(\nu)}(h,q)=E^Q\bigl[(A^{(\nu)}_h-q)\bigr], $$
the {\it normalized time--$t$ price\/} of the Asian option.
To explain the notation, $A^{(\nu)}$ is Yor's twohundred percent 
volatility accumulation process
$$
A^{(\nu)}_h
=\int_0^h e^{2(B_w+\nu w)}dw \, , 
$$
and the normalized parameters are as follows:
$$
\nu= {2\varpi\over \sigma^2}-1,\qquad 
h={\sigma ^2\over 4}(T-t),\qquad 
q=kh\!+\! q^*,
$$
where
$$
k={K\over S_t}, \qquad q^*=q^*(t)=
{\sigma^2\over 4S_t}\left( K\cdot(t\!-\!t_0)
-\int_{t_0}^t S_u\, du \right).
$$
\smallskip
To interpret these quantities,  
$\nu$ is the {\it normalized adjusted interest rate\/},  $h$ is
the {\it normalized time to maturity\/}, which is 
non--negative, and $q$ is the {\it normalized strike price\/}.   
\goodbreak\bigskip
{\bf 5.\quad A reduction of the valuation problem:}\quad 
Computing the normalized time--$t$ price 
%$$ 
%C^{(\nu)}:=
%E^{Q}\big[ \big( A^{(\nu)}_h-q\big)^{\!+}\big]
%$$
of the Asian option reduces to the case where the normalized strike 
price $q$ is positive. 
Indeed, if $q$ is non--positive, Asian options loose their option 
feature, and their normalized time--$t$ price is given by
$$ C^{(\nu)}
:=E^{Q}\big[ \big( A^{(\nu)}_h-q\big)^{\!+}\big]
=E^{Q}\Big[ A^{(\nu)}_h\Big]-q\, . 
$$ 
On applying Fubini's theorem, this last expectation is computed as follows: 
$$
E^{Q}\Big[ A^{(\nu)}_h\Big]
= {e\vbox to 7pt{}^{2h(\nu+1)}-1\over 2(\nu\!+\!1)}
\, ,
$$
for $\nu$ any real number. The right hand side is analytic in $\nu$ 
with its value at $\nu=-1$ being equal to $h$. For computing 
$C^{\raise-2pt\hbox{$\scriptstyle (\nu)$}}$ if $q$ is positive we discuss
two Laguerre series. They in turn crucially depend 
on an  extension of a basic result of Yor's about the Asia density to 
be discussed next. 
\goodbreak\bigskip
{\bf 6.\quad Yor's Asian option density revisited:}\quad This section 
discusses growth and integrability properties of the density of the 
accumulation processes $A^{\raise-1.5pt\hbox{$\scriptstyle (\nu)$}}$. 
This is based on Yor's integral representations of [{\bf Y}, (6.e), p.528]
recalled below and extends his measurability results. Recall 
the Girsanov identity for the normalized time--$t$ price
$$
C^{(\nu)}=E^Q\big[ (A_h^{(\nu)}-q)^+\big]
         = e\vbox to 8pt{}^{-{\scriptstyle \nu^2 h\over\scriptstyle 2}}
          E\big[ (A_h^{(0)}-q)^+ e\vbox to 8pt{}^{\nu W_h}\big]
$$
which follows on making $W_t=\nu t\!+\!B_t$ a Brownian motion and 
dropping reference to the resulting measure. It reduces computing
$C^{\raise-1.5pt\hbox{$\scriptstyle (\nu)$}}$ to a corresponding 
problem in the zero drift accumulation process  
$A^{\raise-1.5pt\hbox{$\scriptstyle (0)$}}$
at the cost of introducing a second stochastic factor in the integrand. 
Define the density $\alpha_{\nu,h}$ on the positive real line by
$$
\alpha_{\nu,h}(x)= e\vbox to 8pt{}^{-{\scriptstyle x^2 \over\scriptstyle 2}}
     \int_0^\infty g_\nu(y) 
             e\vbox to 8pt{}^{-{\scriptstyle xy^2 \over\scriptstyle 2}}
             \psi_{xy}(h)\, dy\, , $$
for any $x>0$. Here $g_\nu$ is the map on the positive real line that
sends any $y>0$ to its $\nu$--th power, and $\psi_\xi(h)$ is for any
positive real number $\xi$ given by:
$$
\psi_\xi(h)=\int_0^\infty 
   e\vbox to 8pt{}^{-{\scriptstyle w^2 \over\scriptstyle 2h}-\xi\cosh(w)}
   \sinh(w)\sin\Big( {\pi w\over h}\Big) dw\, . $$
To  address integrability 
%for Yor's representation, 
%triple integral of [{\bf Y}, (6.e), p.528], 
let $f$ be any continuous function 
on the positive real line, and for any real numbers $a$, $b$  define 
the weighted maps $f_{a,b}$ by $f_{a,b}(x)=\exp(x^{-1}a)x^bf(x)$, 
for any positive real number $x$. Then we have the 
\bigskip
{\bf Proposition:}\quad{\it If $a<1/2$ and $f$ behaves like a power 
map near the origin and towards infinity, the map sending any positive 
real number $x$ to $f(x^{-1})\exp(ax)x^{-b}\alpha_{\nu,h}(x)$ extends 
to an integrable map on the non--negative real line, and we have the 
identity
$$
E\big[ f_{a,b}(A_h^{(0)})g_\nu(e^{W_h})\big]
=c_h\int_0^\infty 
        f_{a,b}\Big( {1\over x}\Big)      
                   %f\Big( {1\over x}\Big) e\vbox to 8pt{}^{ax} x^{-b} 
        \alpha_{\nu,h}(x)\, dx
$$ 
of finite integrals, where $c_h=(2\pi^2h)^{-1/2}\exp((2h)^{-1}\pi^2)$.
Moreover, the map sending any positive real number 
$x$ to $\exp(ax)x^{-b}\alpha_{\nu,h}(x)$ extends to an integrable
and bounded map on the non--negative real line.}
\medskip
{\bf Remark:}\quad Examples of functions $f$ which meet the 
conditions of the Proposition are furnished by any rational 
function with no poles on the positive real line, and by 
any function sending $x\ge 0$ to $(x\!-\!q)^+$ or $(q\!-\!x)^+$.
\medskip
The proof of the Proposition is based on
Yor's triple integral [{\bf Y}, (6.e), p.528] and  his observation 
[{\bf Y}, (6.g), p.529]. Recall the latter implies that for any
non--negative integer $k$ the map $\xi\mapsto \psi_\xi(h)$ is of 
order big Oh of $\xi^k$ with going $\xi$ to zero, while the former
asserts
$$
E\big[ h(A_h^{(0)})g(e^{W_h})\big]
=c_h\int_0^\infty \int_0^\infty h\Big( {1\over x}\Big) 
           e\vbox to 8pt{}^{-{\scriptstyle x\over\scriptstyle 2}} 
g(y) e\vbox to 8pt{}^{-{\scriptstyle xy^2\over\scriptstyle 2}}
\psi_{xy}(h)\, dx\, dy\,  $$ 
for any Borel measurable functions $g$, $h$ on the non--negative real
line into itself.  
In the Yor triple integral for the Proposition you would like to interchange
the order of integration. With the functions $f$ and $\xi\mapsto \psi_\xi(h)$  
possibly negative, Tonelli's theorem does not apply for this. However
it does apply on taking absolute values of the integrands. And if the 
resulting integral can be proved finite this will give integrability 
and justify interchanging the order of integration in the original 
Yor's triple integral using Fubini's theorem. Abbreviating $\gamma=\nu\!-\!b$
and $\delta=1/2\!-\!a$ we are so reduced to prove finiteness of 
$$
J
= 
\int_0^\infty  |f|\Big( {1\over x}\Big) 
           e\vbox to 8pt{}^{-\delta x} D_{0,\infty}(x)\, dx\, , $$
setting 
$$
D_{A,B}(x)= \int_A^B
{ e\vbox to 6pt{}^{-{\scriptstyle \xi^2\over\scriptstyle 2x}}
\over x\vbox to 6pt{}^{\gamma+1}}\cdot \xi^\gamma\cdot |\psi_\xi(h)|\, d\xi\, .
$$
Decompose $J$ into four parts by breaking the two integrations at
the point $1$. The idea is to first majorize the respective inner
integrals $D_{A,B}$ using Yor's above observation, and then use the 
exponential decay of the integrand towards infinity to bound the 
whole summand. We discuss this for the case where $x$ and $\xi$ both
range from $1$ to infinity. Fixing $x>1$ and changing variables 
$\eta =(2x)^{-1/2}\eta$, we have 
$$
\eqalign{
D_{1,\infty}(x)&={\sqrt{2x}\over x\vbox to 6pt{}^{\gamma+1}}
               \int_{1/\sqrt{2x}}^\infty 
               e\vbox to 8pt{}^{-\eta^2}\cdot (\eta \sqrt{2x\, })^c\cdot 
               |\psi_{\eta\sqrt{2x\, }}(h)|\, d\eta.\cr
%               &\le 
%               {\sqrt{2x}\over x\vbox to 6pt{}^{\gamma+1}}
%               \int_{0}^\infty 
%               e\vbox to 8pt{}^{-\eta^2}\cdot (\eta \sqrt{2x\, })^c\cdot 
%               |\psi_{\eta\sqrt{2x\, }}(h)|\, d\eta\, .\cr
} 
$$
To majorize this last integral notice that $(2x)^{1/2}\eta$ is 
greater than or equal to $1$ by construction. On $[1,\infty)$ the 
map $y\mapsto y^\gamma\exp(-y)$ is bounded by a constant $A_\gamma$ 
depending only on $\gamma$. Use this in the defining integral for 
$\psi$ to obtain   
$$
y\vbox to 7pt{}^\gamma |\psi_y(h)|
\le \int_0^\infty 
     e\vbox to 8pt{}^{ -{\scriptstyle w^2\over\scriptstyle 2h}}
     y\vbox to 7pt{}^\gamma e\vbox to 7pt{}^{-y}\sinh(w)\, dw
\le A_\gamma  \int_0^\infty 
     e\vbox to 8pt{}^{ -{\scriptstyle w^2\over\scriptstyle 2h}}
     \sinh(w)\, dw\, . $$ 
Thus there is a constant $B_{\gamma,h}$ depending only on $\gamma$ and $h$
such that $D_{1,\infty}(x)$ sits below
$B_{\gamma,h}\cdot x^{-(\gamma+1/2)}$
for any $x$ in $[1,\infty)$. For the respective summand of $J$ we
so obtain,
$$
\int_1^\infty 
           |f|\Big( {1\over x}\Big) 
           e\vbox to 8pt{}^{-\delta x} D_{1,\infty}(x)\, dx
\le B_{\gamma, h}
          \int_1^\infty 
           |f|\Big( {1\over x}\Big) 
      \cdot {e\vbox to 8pt{}^{-{\scriptstyle\delta\over\scriptstyle 2} x}
            \over x\vbox to 6pt{}^{\scriptstyle\gamma +1/2}}
      \cdot e\vbox to 8pt{}^{-{\scriptstyle \delta\over\scriptstyle 2} x}
\, dx\, .$$
Now apply the above reasoning to the middle factor of the last integral's
integrand and bound it by $(\delta/2)^{\gamma +1/2} A_{-(\gamma +1/2)}$.
On substitution,
$$\eqalign{
\int_1^\infty 
           |f|\Big( {1\over x}\Big) 
           e\vbox to 8pt{}^{-\delta x} &D_{1,\infty}(x)\, dx
\cr &
\le B_{\gamma, h}\cdot \Big({\delta\over 2}\Big)^{\gamma +1/2} 
A_{-(\gamma +1/2)}
  \int_1^\infty 
           |f|\Big( {1\over x}\Big)
\cdot e\vbox to 8pt{}^{-{\scriptstyle \delta\over\scriptstyle 2} x}
\, dx\, ,\cr }$$
and this last integral is finite since $\delta $ is positive by 
construction and $f$ behaves like a power map near the origin by 
hypothesis. The remaining three summands of $J$ can be majorized 
in a similar way, and this completes the proof. 
\goodbreak\bigskip
{\bf 7.\quad First Laguerre series for the Asian option:}\quad 
As a first application of our study of the growth behaviour of 
Yor's Asia density we illustrate a general priciple 
for getting  series expansions of option prices in terms of 
higher moments of the option's control variable. The idea is 
to Laguerre expand in the risk neutral expectation that gives 
this value the taking--the--non--negative--part function, and 
then try to put the expectation through the series so obtained. 
The krux of this {\it Ansatz\/} is to justify the last operation, 
and if such a justification can be found it will be specific to 
the option considered. For the Asian option we have been able to 
give such a justification using the properties of the Asia density 
established in \S 6~Proposition in an essential way. The 
precise result is the following series in terms of the moments 
of the reciprocal of Yor's accumulation process.
\bigskip
{\bf Lemma:}\quad {\it Let $c$ be any positive real number with
$0<qc<1/2$. Setting $X=A_h^{(\nu)}$, the normalized 
time--$t$ price $C^{\raise-1.5pt\hbox{$\scriptstyle (\nu)$}}$ 
of the Asian option has the absolutely convergent series  
$$
C^{(\nu)}
=
{1\over c} \sum_{n=0}^\infty 
c_n\cdot E^Q\Big[X\cdot L_n^\alpha\Big( % A_h^{(\nu)}  L_n^\alpha \Big( 
{qc\over X}\Big) %A_h^{\raise-1.5pt\hbox{$\scriptstyle (\nu)$}} }\Big)
\Big] 
$$
whose any $n$--th coefficient $c_n$ is given by  
$$
c_n=\sum_{k=0}^n {(-1)^{k+1}\over \Gamma(k\!+\!\alpha\!+\!1)}{n\choose k}
\Big( (k\!+\!\alpha)\cdot \gamma(k\!+\!\alpha\!+\!1,c)
+c\vbox to 7pt{}^{ k+\alpha+1} e\vbox to 7pt{}^{-c}\Big),
$$
with $\gamma(s,a)=\int_0^a \exp(-x)x^{s-1}\, dx$ for any $a$, $s>0$ 
the incomplete gamma function, and} 
$$
E^Q\Big[ X\cdot L_n^\alpha \Big(  %A_h^{(\nu)}  L_n^\alpha \Big( 
{qc\over X } \Big) %A_h^{\raise-1.5pt\hbox{$\scriptstyle (\nu)$}} }\Big)
\Big]
=\sum_{k=0}^n (-1)^k 
{\Gamma(n\!+\!\alpha\! +\!1)\over\Gamma(k\!+\!\alpha\! +\!1)}
{(qc)^k\over k!(n\!-\!k)!}
E^Q \Big[ {1\over X^{k-1}} \Big] \, . %\big(A_h^{(\nu)}\big)^{1-k}\Big]\, . 
$$    
\smallskip
{\bf Remark:}\quad In [{\bf D}, Theorem 5.1, p.417f] a Laguerre 
series for the 
the density of the reciprocal of Yor's accumulation process is
given. Use it to compute the expectation defining 
$C^{\raise-1.5pt\hbox{$\scriptstyle (\nu)$}}$ by a formal
term by term integration. The series so obtained is then seen 
to be identical with that of the Lemma.  
\medskip
The proof of the Lemma illustrates the basic  principles
for working with orthogonal polynomials. Take any function $\phi$
on the positive real line which has an $\alpha$--Laguerre series,
so let $w^{1/2}\phi$ be square integrable where $w$ is the weight
$w(x)=\exp(-x)x^\alpha $. From 
[{\bf AAR}, Theorem 6.5.3, p.308] we have mean convergence with 
respect to $w$ of this series to $\phi$. 
Writing $c_k$ for any $k$--th Laguerre coefficient and 
$R_n =\phi-(c_0L_0^\alpha+\cdots +c_nL_n^\alpha)$ for any $n$--th 
order remainder term of the Laguerre series, this is equivalent to  
$$\lim_{n\rightarrow \infty} \int_0^\infty w(x)\cdot
[ R_n(x)]^2\, dx =0\, . $$
The idea for the proving the Lemma is to transcribe the ordinary
convergence of the series to be proved into mean convergence 
of the Laguerre series for $\phi_c(x)=(c\!-\!x)^+$ and apply this
mean convergence result. Thus let $R_n$ now denote any $n$--th 
order remainder term of the Laguerre series for $\phi_c$
and set $\rho_n(x)=xR_n(x^{-1}qc)$. Observing 
$(x\!-\!q)^+=c^{-1}x\phi_c(x^{-1}qc)$, %we have the tautology
$$
E^Q\Big[ \big(X\!-\!q\big)^+\Big]   %A^{(\nu)}_h\!-\!q\big)^+\Big]
-{1\over c}\sum\nolimits_{k=0}^n 
c_k E^Q\Big[ X\cdot L_k^\alpha\Big({qc\over X}\Big)\Big]
%
%A^{(\nu)}_hL_k^\alpha\Big({qc\over A^{(\nu)}_h}\Big)\Big]
=
{1\over c} E^Q\Big[\rho_n\big(X\big)\Big]\, .%A_h^{(\nu)}\big)\Big]\, . 
$$
The proof of the Lemma  so reduces to show that this 
last expectation goes to zero with $n$ to infinity. Using 
the Asian density $\alpha_{\nu,h}$ of \S 6~Proposition, we
have 
$$
E^Q\Big[\rho_n\big(X \big)\Big]   %        A_h^{(\nu)}\big)\Big]
\le c_0 \int_0^\infty \alpha_{\nu,h}(x)\rho_n\Big({1\over x}\Big)\, dx\, ,
$$
with $c_0=c_h\exp(-\nu^2h/2)$. Applying Cauchy--Schwarz,
$$
\eqalign{
\Big|E^Q\Big[\rho_n\big( X\big)\Big]\Big|^2  %A_h^{(\nu)}\big)\Big]\Big|^2
&\le 
c_0^2 \int_0^\infty {1\over x^2} \alpha_{\nu,h}(x)\, dx
      \cdot \int_0^\infty |R_n(qcx)|^2\alpha_{\nu,h}(x)\, dx\cr 
&\le {c_0\over qc} E^Q\big[X^2\big]\cdot %\big(  A_h^{(\nu)}\big)^2\big]\cdot 
\int_0^\infty
\big|R_n(\xi)\big|^2\alpha_{\nu,h}\Big({\xi\over qc}\Big)\, d\xi \, , \cr 
}
$$
on changing variables $\xi=qcx$ for the last equality. Since 
all higher moments of any Yor's accumulation process are finite
using [{\bf Y}, \S 4], we are further reduced to show that
the second integral factor goes to zero with $n$ to infinity. 
This is where the growth properties of the Asian option density 
enter in an essential way. Indeed, using the hypothesis $0<qc<1/2$ 
choose any $qc<\beta<1/2$. Then \S 6~Proposition applies to show 
that $x\mapsto \exp(\beta x)x^{-\alpha} \alpha_{\nu,h}(x)$ is bounded 
on the positive real line, whence
$$\xi\mapsto 
e\vbox to 8pt{}^{\xi}\cdot \Big({\xi \over qc}\Big)^{-\alpha}
\cdot \alpha_{\nu,h}\Big({\xi\over qc}\Big)
$$ 
is bounded on the non--negative real line as well. As a consequence
there is a positive constant $D$ such that 
$$ 
\int_0^\infty
|R_n(\xi)|^2\alpha_{\nu,h}\Big({\xi\over qc}\Big)\, d\xi
\le 
D \int_0^\infty e^{-\xi}\xi^\alpha  |R_n(\xi)|^2 \, d\xi
=
D \int_0^\infty w(\xi)|R_n(\xi)|^2 \, d\xi 
\, .
$$
Now the mean convergence result recalled at the beginning
applies to give that this last integral goes to zero
with $n$ to infinity. Absolute convergence
of the series follows using Riemann's criterion for conditional
convergence: with the above series also any rearragement
of it converges to the same function. 
%$C^{\raise-2pt\hbox{$(\scriptstyle \nu)$}}$. 
This completes the proof of the Lemma.  
\goodbreak\bigskip
{\bf 8.\quad Ladder height density Laguerre series for the Asian 
option:}\quad  A second Laguerre expansion of the value of the Asian 
option has been proposed in [{\bf D}]. The idea is to construct this 
value as a probability density in the strike price, Laguerre expand, 
and try to re--interpret the Laguerre coefficients in terms of the 
reciprocal of Yor's accumulation process. The key notion is that
of the {\it ladder height density\/} $\ol g_Z$ associated to a
probability density function $g_Z$ of any positive $Q$--integrable 
random variable $Z$. It is defined for any $c>0$ by
$$  
\ol g_Z (c)= {1\over E^Q[Z]}\int_c^\infty w\cdot g_Z(w)\, dw\, . $$
%and we denote by $Z_{\ol g}$ the corresponding random variable with 
%$Q$--probability function $\ol g_Z$. 
%\medskip
We have slightly modified the argument of [{\bf D}, \S 7] how this 
applies to Asian options. Using \S 5 assume the modified 
strike price $q$ positive. For any positive real number $c$,
define the $Q$--probability density function $g$ of a positive 
and $Q$--integrable random~variable~by
$$
g(y)={1\over E[Y^{-1}]}\cdot{f_Y(y)\over y}
\qquad\hbox{where}\qquad Y={qc\over A_h^{(\nu)}}\, . 
$$
Then we have the representation of the normalized time--$t$ price 
$C^{\raise-2pt\hbox{$(\scriptstyle \nu)$}}$ of the Asian option
in terms of the ladder height density $\ol{\ol{g}}$ 
associated to $\ol g $ 
$$
C^{(\nu)}=E^Q\big[X\big]-q  %        A_h^{(\nu)}\big]-q
+{q^2c\over 2}E^Q\big[X^{-1}\big]\cdot \ol{\ol{g}}(c) 
%\big(A_h^{(\nu)}\big)^{-1}\big]\cdot \ol{\ol{g}}(c)
%\qquad\hbox{where}\qquad X=A_h^{(\nu)}\, . 
$$    
setting $X=A_h^{(\nu)}$. 
Postponing the computation of $E[X^{-1}]$ to \S 12~Lemma and with $E^Q[X]$ 
computed in \S 5, we are so reduced to computing $\ol{\ol{g}}(c)$. 
The idea is to try to use a Laguerre series for this. Thus consider 
for any positive real number $c$ and any real numbers $\beta$ and 
$\delta$ the formal Laguerre series expansion  
$$
c\vbox to 8pt{}^\beta e\vbox to 8pt{}^{-\delta c}\cdot \ol{\ol{g}}(c)
=
\sum\nolimits _{n=0}^\infty c_n L_n^\alpha(c)$$
whose any $n$--th Laguerre series coefficient is given by
$$
c_n= \sum\nolimits_{k=0}^n {(-1)^k\over \Gamma(k\!+\!\alpha\!+\!1)}
{n\choose k}\cdot I_{\beta,\delta,k}\, . 
$$
where $I_{\beta,\delta,k}$ abbreviates the integral
$$
I_{\beta,\delta,k}=\int_0^\infty y\vbox to 8pt{}^{\alpha+\beta+k}
e\vbox to 8pt{}^{ -(1+\delta)y}\cdot \ol{\ol{g}}(y)\, dy\, .
$$ 
The basic question is under which conditions on $\beta$ and $\delta$ this 
series converges, and this is again crucially based on the analysis of
Yor's density in \S 6. The second question is how to compute the coefficient 
integrals $I_{\beta,\delta,k}$ for these admissible values of $\beta$ 
and $\delta$.  Here is  the following extension of 
[{\bf D}, Theorem 7.1, p.422].
\bigskip
{\bf Proposition:}\quad{\it Assume $q>0$, abbreviate $X=A^{(\nu)}_h$,  
and let  $\beta$ and $\delta$ be such that $\alpha\!+\!2\beta>-1$ and 
$2\delta>-(1\!+\!(2qc)^{-1})$. Then the above $\alpha$--Laguerre series 
converges. If moreover $qc<1/2$, we can choose $\delta=-1$ as an admissible 
such value, and 
$$ I_{\beta,-1,k}={2\over E^Q[X^{-1}]} 
\cdot
{(qc)^{\alpha+\beta+k}
\over (\alpha\!+\!\beta\!+\!k\!+\!1) (\alpha\!+\!\beta\!+\!k\!+\!2)}
E^Q\bigg[ {1\over X^{\alpha+\beta+k+1}}\bigg]\, ,  
$$
for any non--negative integer $k$. 
If $\delta\ne -1$ and $m=\alpha\!+\!\beta\!+\!k$ is any non--negative integer,
$$
\eqalign{
I_{\beta,\delta,k}= 
{2qc\over E^Q[X^{-1}]} {m!\over (1\!+\!\delta)^{m+1}}
&-
{2\over E^Q[X^{-1}]}
\sum_{\l=0}^m
{m!\over (1\!+\!\delta)^{\l+1}}\bigg\{
{E[X^{-1}]\over (1\!+\!\delta)^{m-\l}}
\cr &-\sum_{p=0}^{m-\l}
{(qc)^{m-\l-p}\over (m\!-\!\l\!-\!p)!(1\!+\!\delta)^{p+1}}
E\Big[{1\over X^{m-\l-p-1}}
e\vbox to 9pt{}^{ -{\scriptstyle qc(1+\delta)\over \scriptstyle X}}
\Big]\bigg\}.\cr } 
$$
}
\smallskip
To prove the Proposition, given [{\bf D}, Theorem 5.1, p.418] 
we first have to check square integrability, 
i.e., have to determine which  maps $x\mapsto 
x^{\alpha+2\beta}\exp(-(1+2\delta)x)(\ol{\ol{g}}(x))^2$ are
integrable on $(0,\infty)$. To get a hold on $\ol{\ol{g}}$, 
recall if $X_{\ol g}$ is any random variable 
with $Q$--density function $\ol g $, we have 
$E^Q[X_{\ol g}]\ol{\ol{g}}(x)= Q(\exp(aX_{\ol g})>\exp(ax))$ for any
$a>0$. Hence  Markov's inequality gives 
$E^Q[X_{\ol g}]\ol{\ol{g}}(x)\le \exp(-ax)E^Q[\exp(aX_{\ol g})]$. 
On substitution we are so reduced to characterize finiteness of 
$${E^Q[\exp(aX_{\ol g})]\over E^Q[X_{\ol g}]}\int_0^\infty 
x\vbox to 8pt{}^{\alpha+2\beta}
e\vbox to 8pt{}^{ -(1+2(\delta+a))x}dx\, . 
$$
The integral factor gives the conditions $\alpha\!+\!2\beta>-1$ and 
$1\!+\!2(a\!+\!\delta)>0$. Since $E^Q[X_{\ol g}]=E^Q[X^2]/2$ using 
[{\bf D}, Theorem 2.1, p.411] and all moments of $X$ are finite, we 
have to determine those $a>0$ for which  $E^Q[\exp(aX_{\ol g})]$ is 
finite. Since on unravelling definitions 
$$
E^Q[\exp(aX_{\ol g})]
=
{E^Q[X]\over aqc}
\bigg({1\over E^Q[X]} 
E^Q\Big[ Xe\vbox to 9pt{}^{ {\scriptstyle aqc\over\scriptstyle X}}\Big]
-1\bigg),
$$
this again is an application of \S 6~Proposition. The 
expectation will be finite for any positive real number $a$
such that $aqc<1/2$ yielding $2\delta>-(1\!+\!(2qc)^{-1})$.
This is satisfied for $\delta=-1$ if $qc<1/2$, and in this case
the coefficient integrals $I_{\beta,\delta,k}$ have been computed
in [{\bf D}, Theorem 2.1, p.411]. An analogous computation gives 
the coefficient integrals in the other case which is not contained 
in [{\bf D}]. This completes the proof of the Proposition. 
\goodbreak\bigskip
{\bf 9.\quad Review of two results of Dufresne:}\quad As one of his 
insights into the structure of the accumulation process 
$$
A^{(\nu)}_h=\int_0^h e\vbox to 7pt{}^{ 2(\nu w+B_w)}dw\, , $$
Yor describes in [{\bf Y}, (4.k), p.522] how its higher moments
are determined by certain Gauss hypergeometric functions ${}_2F_1$. 
These results have been extended in [{\bf D}, \S 4] to the higher
moments of the reciprocals of these accumulation processes in 
the following way. First we have the integral representation
$$
E^Q\Big[ {1\over A_h^{(\nu)}}\Big]
=2{ e\vbox to 6pt{}^{-{\scriptstyle \nu^2h\over\scriptstyle 2}}\over
\sqrt{2\pi h^3\, }}
\int_0^\infty ye\vbox to 7pt{}^{-{\scriptstyle y^2\over\scriptstyle 2h}}
{\cosh((\nu-1)y)\over \sinh(y)}\, dy\, , $$
for any real number $\nu$. At the base of this result is Yor's
description of the Laplace transform with respect to time of the 
first moment of $A^{\raise-2pt\hbox{$\scriptstyle (\nu)$}}$ using the
confluent hypergeometric function $\Phi$. This description remains 
valid {\it mutatis mutandis\/} for the  expectation of $r$--th powers 
of $A^{\raise-2pt\hbox{$\scriptstyle (\nu)$}}$ where $r>-1$. The 
identity resulting on inversion is in terms of the Gauss hypergeometric 
function and holds for $r>-3/2$ using analytic continuation. 
\medskip
The higher moments of the reciprocal of    
$A^{\raise-2pt\hbox{$\scriptstyle (\nu)$}}$
are determined recursively from its first moment. Considering
any $k$--th such moment 
$$
m_k(h)= E^Q\big[ \big(A_h^{(\nu)}\big)^{-k}\big]
%E\Big[ \Big( {1\over A_h^{(\nu)}}\Big)^k\Big] 
$$
as a function in the time variable $h>0$, we have
$$ m_k= 2(k-(\nu\!+\!1))m_{k-1}-{1\over k-1}m'_{k-1}\, , $$
for any $k\ge 2$. This is proved using the It\^o Lemma on applying
time reversal. Different angles on this are discussed in [{\bf Y1}],
[{\bf Y2}], and [{\bf Y3}] in particular. 
\goodbreak
\bigskip
{\bf 10.\quad Higher moments as Theta integrals in the basic case:}\quad
It is Theta functions as discussed in \S 3 which provide the 
key for understanding the higher moments of the reciprocal of Yor's 
accumulation process. With higher moments of Yor's accumulation 
processes themselves being determined by hyperbolic sines 
[{\bf Y},4.g, p.519], 
this appearance of theta functions seems unexpected at least.
This section considers the {\it basic case\/} $|\nu-1|\le1$ in which
the higher moments of the reciprocal of any time--$h$ value of 
$A^{\raise-2pt\hbox{$\scriptstyle (\nu)$}}$ 
$$
m_k(h)=  E^Q\big[ \big(A_h^{(\nu)}\big)^{-k}\big]
%E\Big[ \Big( {1\over A_h^{(\nu)}}\Big)^k\Big] 
$$ 
are completely determined by theta functions. Indeed, computing 
them reduces to computing the Theta integrals $\Theta_k$ given by
$$
\Theta_k(h)={1\over 2\sqrt{2\, }}
\int_0^\infty \vartheta\Big({\nu\over 2}\Big| w\Big)
{w^{k-1}\over (wh+1/2)^{k+1/2}}\, dw\, , 
$$
for any $h>0$ and for any positive integer $k$. The precise result
is the following
\bigskip
{\bf Lemma:}\quad {\it Suppose $|\nu-1|\le 1$. For any positive 
integer $n$, we have 
$$
e\vbox to 7pt{}^{{\scriptstyle \nu^2h\over\scriptstyle  2}}
m_n(h)=
a_{n,1}\Theta_1(h)+\cdots +a_{n,n}\Theta_n(h)\, , $$
for any $h>0$, where the coefficients $a_{n,k}$ are recursively determined by
$a_{1,1}=1$ and the recurrence relations
$$
\eqalign{
a_{n+1,1}&= \Big( 2(n\!-\!\nu)+{\nu^2\over 2n}\Big)a_{n,1}\, , \cr 
a_{n+1,n+1}&= \Big( 1\!+{1\over 2n}\big)a_{n,n}\, ,            \cr 
a_{n+1,k}&= \Big( 2(n\!-\!\nu)+{\nu^2\over 2n}\Big)a_{n,k}
+\Big( {k\!-\!1\over n}+{1\over 2n}\Big)a_{n,k-1}\, ,              \cr
}
$$
for any $k$ between $2$ and $n$.}
\medskip
{\bf Remark:}\quad Similar but different charactarizations of
these higher moments have been obtained by Yor and coworkers 
if $\nu=0$, see in particular [{\bf Y1}], [{\bf Y2}], and 
[{\bf Y3}, \S 4].
\medskip
The key insight for the proof is that in Dufresne's integral for 
the expectation recalled in \S 9 the hyperbolic quotient factor of 
the integrand is the Laplace transform at $y^2$ of a theta function
if $|\nu-1|\le 1$. Indeed, for any positive real number $y$ we then
have 
$$ 
{\cosh((\nu\!-\!1)y)\over y\sinh(y)}
=
\int_0^\infty e^{-y^2w}\vartheta\Big( {\nu\over 2}\Big| w\Big)\, dw\, . 
$$
Applying Fubini's theorem the expression of the Lemma for $m_1(h)$ 
follows from this. An induction using the recursion relation for 
$m_k$ of \S 9 then completes the proof.
\bigskip
{\bf 11.\quad Computing the Theta integrals:}\quad This section discusses
two series for computing the Theta integrals $\Theta_n$. From \S 10 
these are recalled to be given by
$$
\Theta_n(h)={1\over 2\sqrt{2\, }}
\int_0^\infty \vartheta\Big({\nu\over 2}\Big| w\Big)
{w^{n-1}\over (wh+1/2)^{n+1/2}}\, dw\, , 
$$
for any $h>0$ and for any positive integer $n$. 
\medskip
{\bf First series:}\quad First we have the 
following absolutely convergent series in terms of the confluent
hypergeometric function of the second kind $\Psi$
$$
\Theta_n(h) 
=
{\Gamma(n)\over 2h^n}\Psi\Big(n,{1\over 2};0\Big)
+
{\Gamma(n)\over h^n}\sum\nolimits_{m=1}^\infty \cos(\pi m \nu)\cdot 
\Psi\Big(n,{1\over 2};{(\pi m)^2\over 2h}\Big)
\, .
$$
It is obtained by term by term integration of the series 
which results from applying the Jacobi transformation formula of \S 3 
to the Theta function factor in $\Theta_n(h)$. However, it is 
problematic in that confluent hypergeometric functions 
are most complicated functions and their implementations are often 
flawed. Similar series for $\nu=0$ are in [{\bf Y3}, \S4].  
\medskip
{\bf Second series:}\quad Still there is a second series in 
terms of less problematic functions. It is the absolutely convergent 
series
$$
\Theta_n(h) 
=
\sum_{m=0}^\infty c_{B,m}(h)+ \sum_{m\in {\bf Z}} d_{B,m}(h)
$$
whose coefficients are for any positive real number $B$ given 
as follows: 
$$
\displaylines{
c_{B,0}(h)
=
{1\over 2h^n}\sum_{k=0}^{n-1} 
{(-1)^{n-1-k}\over n\!-\!k\!-{\displaystyle 1\over\displaystyle 2}}
 {n\!-\!1\choose k} {1\over (2hB\!+\!1)^{n-k-1/2}},\cr 
\noalign{\vskip5pt}
c_{B,m}(h)
=
\gamma_m
\sum_{k=0}^{n-1} (-2)^k {n\!-\!1\choose k} {h^k\over (\pi m)^{2k}}
\Big[ C^{(1)}_{m,k}+ C^{(2)}_{m,k}\Big]
\cr}
$$
for any positive intgeger $m$, abbreviating
$$
\displaylines{
\gamma_m=
{\cos(\pi m \nu)\over \sqrt{2\, }}{(-1)^{n-1}\over 2^{n-1}}
{(\pi m )^{2n-1}\over h^{2n-1/2}}\, , \cr 
\noalign{\vskip4pt}
 C^{(1)}_{m,k}
=
(-1)^{n-k-1} \sum_{\l=0}^{n-k-1}
(-1)^\l 
{
%(1/2)_\l\over (1/2)_{n-k}}
\big({\displaystyle 1\over\displaystyle 2}\big)_\l\over 
\big({\displaystyle 1\over\displaystyle 2}\big)_{n-k}}
\cdot
{e\vbox to 8pt{}^{ \scriptstyle -(\pi m)^2B}
\over \Big((\pi m)^2\Big( 
             B+{\displaystyle 1\over\displaystyle 2h}\Big)\Big)^{\l+1/2}
}\,,\cr
\noalign{\vskip4pt}
 C^{(2)}_{m,k}
=
\sqrt{\pi\, } {(-1)^{n-k}\over 
%(1/2)_{n-k}}
\big({\displaystyle 1\over\displaystyle 2}\big)_{n-k}}
%\sqrt{\pi\, } 
\cdot e\vbox to 9pt{}^{ -(\pi m)^2B}
W\Big( \pi m 
%\sqrt{ B+{1\over 2h}\, }\, 
\sqrt{ B+(2h)^{-1}\, }\, 
\Big).\cr
}
$$
%
%{\cos(\pi m \nu)\over \sqrt{2\, }}&{(-1)^{n-1}\over 2^{n-1}}
%{|\pi m |^{2n-1}\over h^{2n-1/2}}
%\sum_{k=0}^{n-1} (-2)^k {n\!-\!1\choose k} {h^k\over (\pi m)^{2k}}
%\bigg\{\cr 
%\noalign{\vskip3pt}
%&
%(-1)^{n-k-1} \sum_{\l=0}^{n-k-1}
%(-1)^\l 
%{
%%(1/2)_\l\over (1/2)_{n-k}}
%\big({\displaystyle 1\over\displaystyle 2}\big)_\l\over 
%\big({\displaystyle 1\over\displaystyle 2}\big)_{n-k}}
%\cdot
%{e\vbox to 8pt{}^{ \scriptstyle -(\pi m)^2B}
%\over \Big((\pi m)^2\Big( 
%             B+{\displaystyle 1\over\displaystyle 2h}\Big)\Big)^{\l+1/2}
%}\cr 
%\noalign{\vskip3pt}
%+&
%\sqrt{\pi\, } {(-1)^{n-k}\over 
%%(1/2)_{n-k}}
%\big({\displaystyle 1\over\displaystyle 2}\big)_{n-k}}
%%\sqrt{\pi\, } 
%\cdot e\vbox to 9pt{}^{ -(\pi m)^2B}
%W\Big( |\pi m| 
%%\sqrt{ B+{1\over 2h}\, }\, 
%\sqrt{ B+(2h)^{-1}\, }\, 
%\Big)\bigg\}\cr 
%}
%$$
%for any positive integer $m$. 
If $m$ is any integer different from $-\nu/2$, we have
$$
%\eqalign{
d_{B,m}(h)=
{(-2)^{n-1}\over \sqrt{\pi\, }}
\sum_{k=0}^{n-1} 
D_{m,k}^{(1)} + D_m^{(2)}
$$
abbreviating
$$
\displaylines{
D_{m,k}^{(1)}
=
(-1)^k {
\big({\displaystyle 1\over\displaystyle 2}\big)_k\over 
\big({\displaystyle 1\over\displaystyle 2}\big)_{n}}
\cdot
{
\Big(m\!+\!{\nu\over 2}\Big)^{2(n-k-1)}\over 
\Big({1\over B}\!+\!2h\Big)^{k+1/2}}
\cdot e\vbox to 10pt{}^{ -{\scriptstyle 1\over\scriptstyle B}
                   \big( m+{\scriptstyle \nu\over\scriptstyle 2}\big)^2}\, 
\cr
\noalign{\vskip4pt}
D_m^{(2)}
=
2^{n-1}{(-1)^{n}\over 
\big({\displaystyle 1\over\displaystyle 2}\big)_{n}}
\Big| m\!+\!{\nu\over 2}\Big|^{2n-1}
e\vbox to 10pt{}^{ -{\scriptstyle 1\over\scriptstyle B}
                   \big( m+{\scriptstyle \nu\over\scriptstyle 2}\big)^2}
W\bigg(
\Big| m\!+\!{\nu\over 2}\Big| \sqrt{B^{-1}\!+\!2h\, }\bigg)
\, , \cr 
}
$$
while for $m=-\nu/2$ any integer, we have
$$
d_{B,-{\scriptstyle \nu\over\scriptstyle 2}}
=
{2^{n-1}\over \sqrt{\pi\, }
\Big(n\!-\!{\displaystyle 1\over\displaystyle 2}\Big)} 
\cdot
{B^{n-1/2}\over (2hB\!+\!1)^{n-1/2}}
\, . $$
Here the function $W$ is given by $W(z)=\exp(z^2)\erfc(z)$, 
for any complex number $z$, and $(\lambda)_k$ is the Pochhammer symbol
of any complex number $\lambda$ recursively given by $(\lambda)_0=1$ and 
$(\lambda)_{k+1}=(\lambda\!-\!k)\cdot (\lambda)_k$ for any non--negative 
integer $k$. The idea for this series is to optimize 
the convergence behaviour using the Jacobi transformation formula of 
\S 3. It is thus obtained by breaking the integral defining 
$\Theta_n(h)$ at $B$ and integrating over $[B,\infty)$ now the 
Jacobi transform of the Theta function in $\Theta_n(h)$.   
\goodbreak
\bigskip
{\bf 12.\quad Higher moments in the general case:}\quad Computing 
the moments of the reciprocal of any time--$h$ value of 
$A^{\raise-2pt\hbox{$\scriptstyle (\nu)$}}$ 
$$
m_k(h)=  E\big[ \big(A_h^{(\nu)}\big)^{-k}\big]
$$
for arbitrary indices $\nu$ is not solely in terms of the Theta 
integrals $\Theta_k(h)$ of \S 11. Additional  non--Theta correction 
terms are required. The precise result is the following
\bigskip
{\bf Lemma:}\quad {\it Let $n^*$ be the smallest integer $m$ such that
$2m\!+\!1\!-|\nu\!-\!1|$ is non--negative. For any positive 
integer $N$, we then have 
$$\eqalign{
e\vbox to 7pt{}^{{\scriptstyle \nu^2h\over\scriptstyle  2}}
m_N(h)=&
\sum\nolimits_{k=1}^{N} a_{N,k}\Theta_k(h)\cr 
&+\sum\nolimits_{n=0}^{n^*-1} \sum\nolimits_{k=0}^{N-1}
\Big\{ b_{N,k}(h)\cdot C_{n,k}(h)- a_{N,k+1}\cdot D_{n,k+1}(h)\Big\}
\, , \cr }$$
for any $h>0$.}
\medskip 
{\bf The functions $C_{n,k}$ and $D_{n,k}$:}\quad 
The functions $C_{n,k}$ and $D_{n,k}$ are for any positive
real number $h$ given as follows. First, 
$$
C_{n,k}(h)
=
h\vbox to 9pt{}^{{\scriptstyle 1\over\scriptstyle 2}} 
\bigg({\beta_n h\over \sqrt{2\, }}\bigg)^{2k+1}\,  
\sum_{\l=0}^{2k+1} (-1)^{\l-1}
{ 2k\!+\!1\choose \l}\cdot
W_\l\bigg(  {\beta_n\sqrt{2h\,}\over 2}\bigg)\cdot
\bigg( {\beta_n \sqrt{2h\, }\over 2}\bigg)^{-\l}
%W_\l\bigg(  {\beta_n\sqrt{2h\,}\over 2}\bigg)\, . 
$$
Here $W_\l$ is given in terms of the complementary incomplete 
gamma function $\Gamma(a,x)$ by 
$W_\l(x)=\exp(x^2)\cdot \Gamma((\l\!+\!1)/2, x^2)$
for any real number $x$, and $\beta_n=2n\!+\!1\!-|\nu\!-\!1|$. 
Then, 
$$
D_{n, k}(h)
=2\vbox to 9pt{}^{-1}
{(-2\gamma_n)^{k}
%2^{\scriptstyle k+{\scriptscriptstyle 1\over\scriptscriptstyle 2}}
\over 
\big({\displaystyle 1\over\displaystyle 2}\big)_{k}}
\Bigg[{1\over\sqrt{\gamma_n\, }} 
%\over \big({\displaystyle 1\over\displaystyle 2}\big)_{k}}
%\cdot 
%\gamma_n^{k-{1\over2}}
W\Big(\sqrt{ 2h\gamma_n\, }\Big)
-{\sqrt{h\, }\over \sqrt{ \pi\, }}
\sum_{\l=0}^{k-1} 
\big({\displaystyle 1\over\displaystyle 2}\big)_\l
%\over \big({\displaystyle 1\over\displaystyle 2}\big)_{k}}
\cdot 
{(-1)^\l \over 2\cdot (h\gamma_n)^{\l+1}}
\Bigg]
$$
where $W$ is the function $W_0$, abbreviating 
$\gamma_n=4^{-1}(|\nu\!-\!1|\!-\!1\!-2n)^2$, and with 
$(\lambda)_k$ the Pochhammer symbol of any complex number $\lambda$ 
recursively given by $(\lambda)_0=1$ and 
$(\lambda)_{k+1}=(\lambda\!-\!k)\cdot (\lambda)_k$ for any 
non--negative integer $k$.
\medskip
{\bf The coefficients $a_{n,k}$ and $b_{n,k}$:}\quad 
The coefficients $a_{n,k}$ are those of \S 10~Lemma. The 
coefficient functions  $b_{n,k}(h)$ have the form
\belowdisplayskip=0pt
$$
\displaylines{
b_{n,n-1}(h)= b_{n,n-1,n-1} h\vbox to 9pt{}^{-( {3\over2}+2(n-1))}\cr 
\noalign{and for any non--negative integer $k$ less than or equal to 
$n\!-\!2$}
b_{n,k}(h)
=
\sum\nolimits_{\l=0}^{n-1}
b_{n,k,\l}h\vbox to 9pt{}^{-( {3\over2}+k+\l)}
.\cr}
$$
%for $k$ from zero to $n\!-\!2$. 
\goodbreak \belowdisplayskip=8pt plus 5pt minus 3pt
They are recursively defined by
$$ 
\displaylines{
b_{1,0}(h)={2\over \sqrt{\pi\, }} h\vbox to 9pt{}^{-{3\over2}},
\cr
\noalign{\vskip3pt}
b_{2,0}(h)={2\over 2\pi}\Big( 2(1\!-\!\nu)\!+\!{\nu^2\over 2}\Big)
            h\vbox to 9pt{}^{- {3\over2}}
 +{3\over \sqrt{2\pi}}h\vbox to 9pt{}^{- {5\over2}}
\qquad \hbox{and}\qquad
b_{2,1}(h)= -{2\over \sqrt{2\pi}} h\vbox to 9pt{}^{-{7\over2}},\cr
}
$$
and  the following coefficient recurrence rules for $n\ge2$ 
$$
\leqalignno{
&b_{n+1,n,n}=-{1\over n} b_{n,n-1,n-1}\, , &\cr
\noalign{\vskip3pt}
&b_{n+1,n-1,\l}=
%\eqalign{ 
\cases{0 \hfill &\quad $\l=0$\cr
\noalign{\vskip3pt}
          -{\displaystyle 1\over\displaystyle n}
          b_{n,n-2,\l-1}\hfill&\quad  $\l=1,\ldots,n\!-\!2$\cr
\noalign{\vskip3pt}
          \Big(2(n\!-\!\nu)\!+\!
          {\displaystyle \nu^2\over\displaystyle 2n}\Big)b_{n,n-1,n-1}
          -{\displaystyle 1\over\displaystyle n}b_{n,n-2,n-1}\hfill
          & \quad $\l=n\!-\!1$\cr
\noalign{\vskip3pt} 
          {\displaystyle 1\over\displaystyle n}
        \Big( {\displaystyle 3\over\displaystyle 2}\!
         +\!2(n\!-\!1)\Big)b_{n,n-1,n-1}
          -{\displaystyle 1\over\displaystyle n}b_{n,n-2,n-1}\hfill
       & \quad $\l=n$.\cr 
         }&\cr 
}
$$
For $k$ any positive integer less than or equal to $n\!-\!2$ we have 
$$
\leqalignno{
&b_{n+1,k,\l}\cr &=
\cases{ 
    \Big(2(n\!-\!\nu)\!
    +\! {\displaystyle \nu^2\over\displaystyle 2n}\Big)b_{n,k,0}
        &\quad $\l=0$\cr 
\noalign{\vskip3pt}
         \Big(2(n\!-\!\nu)\!
         +\! {\displaystyle \nu^2\over\displaystyle 2n}\Big)b_{n,k,\l}
        +{\displaystyle 1\over\displaystyle  n}
          \Big( {\displaystyle 3\over\displaystyle 2} 
                     \!+\!k\!+\!\l\!-\!1\Big)b_{n,k,\l-1}
        -{\displaystyle 1\over\displaystyle n}
          b_{n,n-1,\l-1}
        &\quad  $\l=1,\ldots,n\!-\!1$\cr
\noalign{\vskip3pt}
         {\displaystyle 1\over\displaystyle n}
        \Big( {\displaystyle 3\over\displaystyle 2}\!
         +\!k\!+n\!-\!1\Big)b_{n,k,n-1}
          -{\displaystyle 1\over\displaystyle n}b_{n,n-1,n-1}\hfill
       & \quad $\l=n$,\cr 
     }&\cr
\noalign{\vskip3pt}
\noalign{and}
\noalign{\vskip3pt}
&b_{n+1,0,\l}=
\cases{ 
    \Big(2(n\!-\!\nu)\!
    +\! {\displaystyle \nu^2\over\displaystyle 2n}\Big)b_{n,0,0}
        &\quad $\l=0$\cr 
\noalign{\vskip3pt}
         \Big(2(n\!-\!\nu)\!
         +\! {\displaystyle \nu^2\over\displaystyle 2n}\Big)b_{n,0,\l}
        +{\displaystyle 1\over\displaystyle  n}
          \Big( {\displaystyle 3\over\displaystyle 2} 
                     \!+\!\l\!-\!1\Big)b_{n,0,\l-1}
        &\quad  $\l=1,\ldots,n\!-\!1$\cr
\noalign{\vskip3pt}
         {\displaystyle 1\over\displaystyle n}
        \Big( {\displaystyle 3\over\displaystyle 2}\!
         +n\!-\!1\Big)b_{n,0,n-1}
       & \quad $\l=n$.\cr 
     }&\cr
}
$$
\goodbreak
\smallskip
{\bf 13.\quad An example:}\quad This section considers the example of 
valuing an Asian option on a 
non--dividend--paying stock with  an interest rate $\varpi=r=9\%$ p.a.,
a volatility $\sigma=30\%$ p.a., and  with time to maturity $1$ year
from today time $0$. Suppose today's price of the stock $S_0$ is equal
to $100$, and the option has been issued today at par, i.e., $K=S_0$. Then 
$ \nu=1.0$, $h=q= 0.0225$,  and we have for the Theta integrals 
$\Theta_n(h)$ of \S 10: 
$$ 
\vcenter{
\vbox{\offinterlineskip 
\hrule
\halign{&\vrule #& \strut \enspace $\hfil#$\enspace\cr
&n &&  n_{c,30}\hfil  && n_{d,30}\hfil
 && \hfill\Theta_n(h)\hfill &\cr  
\noalign{\hrule}
& 1&&4 && 4 &&44.2790749547 &\cr 
\noalign{\hrule}
& 3&&4 && 4 &&46822.1151330265&\cr
\noalign{\hrule}
& 5&&4 && 4 &&70467169.1745463509&\cr
\noalign{\hrule}
& 7&&4 && 4 &&116806496442.2335048967&\cr
\noalign{\hrule}
& 9&&4 && 4 && 202679488380050.1455046663&\cr
\noalign{\hrule}
& 11&&4 &&4 && 361222145817967427.8404500431&\cr
\noalign{\hrule}
& 13&&4 &&4 && 655202294765472642150.1383863999&\cr
\noalign{\hrule}
& 15&&4 &&4  &&1203316829766524026893611.5309617920&\cr
\noalign{\hrule}
& 17&& 4 &&4 &&  2230542763985591766299300393.1676737371&\cr
\noalign{\hrule}
& 19&&4 &&4 && 4164445349148560860975785222305.1797500863&\cr
%\noalign{\hrule}
%& 21&& && && 7819648923423821998234320571395515.9134049386&\cr
\noalign{\hrule}
\noalign{\vskip5pt}
\multispan{9} \hfill \hbox{{\ninebf Table 1}\enspace 
\ninerm  Growth of Theta integrals ${\Theta_n(h)}$.}\hfill\cr 
} 
}}
$$
We have computed these values using the second series of \S 11 with
$B=0.3$, and denote by $n_{c,30}$ and  $n_{d,30}$ those indices 
$m$ for which the first $30$ decimal places of the partial sums 
$c_{B,0}+\cdots +c_{B,m}$ and $d_{B,0}+\cdots +d_{B,m}$ respectively are 
correct. The calculations have been done on a HP Visualize C200 using 
the GP/PARI CALCULATOR Version 2.0.11 (beta) of C. Batut, K. Belabas, 
D. Bernardi, H. Cohen and M. Olivier with 150D precision. Such a 
precision was desired since the coefficients $a_{n,k}$ of \S 10~Lemma
grow like $k$ factorial. Using this result we compute the following 
values of the $n$--th moments $m_n(h)$ of the reciprocal of Yor's 
accumulation process $A^{\raise-2pt\hbox{$\scriptstyle (\nu)$}}$ 
at time $h$:
$$ 
\vcenter{
\vbox{\offinterlineskip 
\hrule
\halign{&\vrule #& \strut \enspace   $\hfil#$\enspace\cr
&n &&   \hfill m_n(h)\hfill &\cr  
\noalign{\hrule}
&1&& 43.7837269185&\cr 
\noalign{\hrule}
&3&& 91741.4898743896&\cr 
\noalign{\hrule}
&5&& 215689033.4795106594&\cr 
\noalign{\hrule}
&7&& 566731384819.7874280500&\cr 
\noalign{\hrule}
&9&& 1657943864684789.3944805424&\cr 
\noalign{\hrule}
&11&& 5380692663910949427.6422561855&\cr 
\noalign{\hrule}
&13&& 19305805471114617878830.4144899460&\cr 
\noalign{\hrule}
&15&& 76329904185193047892084144.9033200670&\cr 
\noalign{\hrule}
&17&& 331513545523453183819373799801.9730752657&\cr 
\noalign{\hrule}
&19&& 1576936784103137595901225118618858.9370371494&\cr 
\noalign{\hrule}
\noalign{\vskip5pt}
\multispan{5} \hfill \hbox{{\ninebf Table 2}\enspace 
\ninerm  Growth of the moments ${m_n}(${\nineit h}).} \hfill\cr 
} 
}}
$$
Given this rapid growth of the negative moments the question is 
if and when convergence of the Laguerre series is rapid enough so that
only the first few of their terms are needed to reproduce the price 
of the Asian option. Principally, convergence of  Laguerre series 
is effected by the oscillation of the Laguerre polynomials $L_n^\alpha$ 
in the range of positive real numbers  smaller than 
$4n\!+2(\alpha\!+\!1)$. In both of our Laguerre series this is controlled by 
the parameter $c$ in the sense that the smaller $c$ the more of
their terms are in their respective oscillatory range, and the better 
should be convergence. The Black--Scholes price of our Asian option
is given as $C^{BS}=8.83$ in [{\bf RS}, Table 1.3] from which its normalized
time--$t$ price 
$C^{\raise-1.5pt\hbox{$\scriptstyle (\nu)$}}=0.002173850758$ 
is obtained on divison by $4061.916379$. In the sequel we consider
the ladder height density series of \S 8~Proposition with $\delta=-1$ 
and $\beta=0$ and take $\alpha=0$, i.e., compute with the classical 
Laguerre polynomials. Using the obvious notation, we record in the 
following examples the values the series give for these two prices 
after summing the first $n\!+\!1$ terms, and give the respective amounts 
$\Delta_n$ that have to be added to obtain their above target values. 
First we choose for $c$ a ``reasonable'' value, to improve the series' 
sprint qualities so to speak.
$$ 
\vcenter{
\vbox{\offinterlineskip 
\hrule
\halign{&\vrule #& \strut \enspace   $\hfil#$\enspace\cr
&n &&   \hfill C^{\raise-1.5pt\hbox{$\scriptstyle (\nu)$}}_n\hfill &&
 \hfill\Delta^{\raise-1.5pt\hbox{$\scriptstyle (\nu)$}}_n \hfill&& 
 \hfill C_n^{BS} \hfill&& \hfill \Delta_n^{BS} \hfill&\cr  
\noalign{\hrule}
&1&&0.0041430388 && 0.0019691881 &&
16.82867729 &&7.9986772963 &\cr 
\noalign{\hrule}
&3 && 0.0029964651&& 0.0008226144
&& 12.17139078 && 3.3413907866 &\cr 
\noalign{\hrule}
&5&& 0.0030265291 && 0.0008526784
&&12.29350849 &&3.4635084905  &\cr 
\noalign{\hrule}
&7&&0.0027607453&&0.0005868946
&& 11.21391670 && 2.3839167088 &\cr 
\noalign{\hrule}
&9&&0.0024684544 &&0.0002946036
&&10.02665531 &&1.1966553198 &\cr 
\noalign{\hrule}
&11&&0.0023864317 &&0.0002125809
&& 9.69348587 &&0.8634858736 &\cr  
\noalign{\hrule}
&13&&0.0023930119 &&0.0002191611
&&9.72021406 &&0.8902140555 &\cr 
\noalign{\hrule}
&15&&0.0023544102 &&0.0001805595
&&9.56341738 &&0.7334173791 &\cr 
\noalign{\hrule}
&17&&0.0022627357 &&0.0000888849
&&9.19104301 &&0.3610430106&\cr 
\noalign{\hrule}
&19&& 0.0021738508 &&4.049\times 10^{-15}
&& 8.83000000 &&1.645\times 10^{-11}&\cr 
\noalign{\hrule}
\noalign{\vskip5pt}
\multispan{11} \hfill \hbox{{\ninebf Table 3}\enspace 
\ninerm  Ladder height series with $c=1.367054258545$}\hfill\cr 
} 
}} %\cr  
$$
As to be expected for this magnitude of $c$, there is a steady 
convergence of the series. One should expect it to converge  
mildly oscillating. However, mind the following fallacy. 
It is not sufficient to have the  the normalized time--$t$ price 
$ C^{\raise-1.5pt\hbox{$\scriptstyle (\nu)$}}$ correct up to the 
the first two or three decimal places only! Indeed,  
we have this accuracy for $ C^{\raise-1.5pt\hbox{$\scriptstyle (\nu)$}}$ 
from the outset while there are still big deviations from the 
Black--Scholes price $C^{BS}$. These phenomenons
are much more pointed if a larger value for $c$ is chosen. Indeed
notice how the following series actually start with almost the 
correct values, then wildly oscillate, and only begin to calm down
after fifteen terms or so.
$$
\vcenter{
\vbox{\offinterlineskip 
\hrule
\halign{&\vrule #& \strut \enspace   $\hfil#$\enspace\cr
&n &&   \hfill C^{\raise-1.5pt\hbox{$\scriptstyle (\nu)$}}_n\hfill &&
 \hfill\Delta^{\raise-1.5pt\hbox{$\scriptstyle (\nu)$}}_n \hfill&& 
 \hfill C_n^{BS} \hfill&& \hfill \Delta_n^{BS} \hfill&\cr  
\noalign{\hrule}
&1&&0.0020571929 &&-0.0001166579
&&8.35614534 &&-0.473854662&\cr 
\noalign{\hrule}
&3&&0.0045737040 &&0.0023998531
&&18.57800280 &&9.748002806&\cr 
\noalign{\hrule}
&5&&0.0044363733 &&0.0022625225
&&18.02017734 &&9.190177342&\cr
\noalign{\hrule}
&7&& -0.0055198391 &&-0.0076936899
&& -22.42112502 &&-31.251125020&\cr
\noalign{\hrule}
&9&&0.0028067490 &&0.0006328989
&&11.40077951 &&2.570779513&\cr
\noalign{\hrule}
&11&&0.0034914141 &&0.0013175634
&&14.18183220 && 5.351832204&\cr
\noalign{\hrule}
&13&&0.0017503210 &&-0.0004235298
&& 7.10965752 &&-1.720342482&\cr
\noalign{\hrule}
&15&&0.0015816389 &&-0.0005922118
&& 6.42448511 &&-2.405514894&\cr
\noalign{\hrule}
&17&&0.0017276113 &&-0.0004462395
&& 7.01741241 &&-1.812587586&\cr
\noalign{\hrule}
&19&&0.0017603640 &&-0.0004134868
&& 7.15045128 &&-1.679548715&\cr
\noalign{\hrule}
\noalign{\vskip5pt}
\multispan{11} \hfill \hbox{{\ninebf Table 4}\enspace 
\ninerm  Ladder height series with $c=6$}\hfill\cr 
} 
}}
$$
In summary, we are now able to efficiently compute the negative 
moments of Yor's accumulation processes that enter into the 
Laguerre series for the Asian option. However, there has still 
work to be done in characterizing optimal choices of the 
convergence parameter $c$. 
\goodbreak\bigskip
{\bf 14.\quad Epilogue:}\quad In this paper we have concentrated
on the first of the two ladder height density Laguerre series of
\S 8~Proposition. Working with the second of these requires to 
have for certain positive real numbers $\alpha$, $\beta$ at least 
computable expressions for the expectations
$$ f_{\alpha,\beta}(t)
=E^Q\Big[ X^{-\alpha} e\vbox to 9pt{}^{ 
-{\scriptstyle \beta\over\scriptstyle 2 X}}\Big]
$$  
for any time $t>0$ setting $X=A_t^{(\nu)}$. Adapting
the argument of [{\bf Y}, \S 3, p.515f], their Laplace transform
with respect to time is computed in terms of modified 
Bessel functions $I_\mu$ and MacDonald functions $K_\mu$ 
%discussed in [{\bf L}, Chapeter 5] 
as follows
$$
\la(f_{\alpha,\beta})(z)
=
{2\over \eta^{\alpha-\nu}}\int_0^\infty \xi^{\alpha-1} 
I_{\sqrt{2z\, }}(\xi)K_{\alpha-\nu}(\eta\xi)\, d\xi
$$
setting $\eta=((1\!+\!\beta)/2)^{1/2}$ and with $z$ any complex number
whose real part is positive and such that $\re((2z)^{1/2})>1\!-\!\nu$.
In analogy to the modified Weber--Schafheitlin integrals of 
[{\bf W}, p.410] these Laplace transforms  can be expressed using the Gauss 
hypergeometric function ${}_2F_1$. On Laplace inversion 
$f_{\alpha,\beta}$ is seen to be given by the~contour~integral
$$
f_{\alpha,\beta}(t)=
{2\over \eta^{\alpha-\nu}} {1\over 2\pi i}\int_{\log C_R} 
{w\over \sqrt{2\pi t^3\, }} 
e\vbox to 9pt{}^{-{\scriptstyle w^2\over \scriptstyle 2t}}
\int_0^\infty \xi^{\alpha-1} K_{\alpha-\nu}(\eta\xi) 
 e\vbox to 9pt{}^{\xi\cosh(w)} d\xi\, dw\, . 
$$
Here $\log C_R$ is the following path of integration 
\input xy
\xyoption{arc}
\xyoption{arrow}
$$
\xy     0;<1mm,0mm>            %%% setzt Koord--System: Einheit ist 1mm...
                \ar(0,0);(80,0),
                \ar(60,-10);(40,-10),
        \ar@{-}(40,-10);(20,-10)*+\dir{*}*+!RU{{\scriptstyle\log R-i\pi}},
        \ar@{-}(20,-10);(20,10)*+\dir{*}*+!RD{{\scriptstyle \log R+i\pi}},
                \ar(20,10);(40,10),
                \ar@{-}(40,10);(60,10),\ar@{},
                (20,0)*+\dir{*}*+!RU{{\scriptstyle\log R}}
\endxy{}
$$
\vskip-.1cm %\parein{10pt}{}
\centerline{\ninerm\hbox to 20pt{} {\ninebf Figure 1}\enspace The 
contour $\hbox{\ninerm log{\nineit C}}_{\scriptscriptstyle R}$}
\vskip.3cm
where the positive real number $R$ has to be chosen so big that the 
real parts of $w^2$ and $w$ of any of its elements $w$ are positive 
and such that $\eta$ is bigger than $\cosh(\log R)$. Under additional 
hypotheses, from [{\bf E}, 7.7.3 (26), p.50] the inner integral can 
be expressed using the Gauss hypergeometric function ${}_2F_1$. 
Alternative angles on these expectations, like their representation as 
Kantorovich--Lebedev transforms, are discussed in [{\bf AG}]. 
\goodbreak
\medskip
\bigskip
{\bf References}
\medskip
{\baselineskip=9.5pt
\ninerm
\parein{35pt}{\ninebf [AAR]} G.E. Andrews, R. Askey, R. Roy: {\nineit 
        Special functions\/}, Encyclopedia of Mathematical Sciences 71, 
        Cambridge UP 1999.
\paraus
\parein{35pt}{\ninebf [AG]} L. Alili, J.--C. Gruet: An explanation 
    of a generalized Bougerol's identity, pp.15--33 in [{\ninebf YE}].
\paraus
\parein{35pt}{\ninebf [D]} D. Dufresne: Laguerre series for Asian and
        other options, {\nineit Math. Finance\/} {\ninebf 10} (2000), 407--28.
\paraus
\parein{35pt}{\ninebf [E]} A. Erd\'elyi et al: {\nineit Higher transcendental
functions II\/}, Krieger reprint 1981.
\paraus
\parein{35pt}{\ninebf [GY]} H. Geman, M. Yor: Bessel processes, Asian  
        options, and perpetuities, {\nineit Math. Finance\/} 
       {\ninebf 3} (1993), 349-375.
\paraus
\parein{35pt}{\ninebf [L]} N.N. Lebedev:  {\nineit Special functions and
        their applications\/}, Dover Publications 1972.
\paraus
\parein{35pt}{\ninebf [M]} D. Mumford, {\nineit TATA lectures on Theta},
Progress in Mathematics 28, 43, 97, Birkh\"auser. 
\paraus
\parein{35pt}{\ninebf [RS]} L.C.G. Rogers, Z. Shi: The value of an Asian 
       option, {\nineit J. Appl. Prob.\/} {\ninebf 32} (1995), 
       1077--88.
\paraus
\parein{35pt}{\ninebf [S]} G. Sansone: {\nineit Orthogonal functions\/},
       Dover Publications 1991.
\paraus
\parein{35pt}{\ninebf[\bf SE]} M. Schr\"oder: On the valuation of
    arithmetic--average Asian options: explicit formulas, Universit\"at
    Mannheim, M\"arz 1999.
\paraus
\parein{35pt}{\ninebf [W]} G.N. Watson: {\nineit A treatise on the 
      theory of Bessel functions\/}, 2nd ed. 1944, reprinted 1995,
       Cambridge UP. 
\paraus
\parein{35pt}{\ninebf [Y]}  M. Yor: On some exponential functionals of 
        Brownian motion, {\nineit Adv. Appl. Prob.} {\ninebf 24}
        (1992), 509--531.
\paraus
\parein{35pt}{\ninebf [YE]}  M. Yor [ed.]: {\nineit Exponential 
    functionals and principal values related to 
    Brownian motion\/}, Revista Mathem\'atica Iberoamericana, 
     Madrid 1997. 

\paraus
\parein{35pt}{\ninebf [Y1]} M. Yor et al: On the positive and negative 
        moments of the integral of geometric Brownian motion, 
        {\nineit Statistics and Probability Letters\/} {\ninebf 49} 
        (2000), 45--52.
\paraus
\parein{35pt}{\ninebf [Y2]} M. Yor et al: On striking identities about the 
        exponential functionals of the Brownian bridge and Brownian motion,
        to appear {\nineit Periodica Math. Hung.} 
\paraus
\parein{35pt}{\ninebf [Y3]} M. Yor et al: The law of geometric Brownian
        motion and its integral, revisited, pre--print Paris VI, 
        October 2000.
\paraus
}

\vfill
\eject
\headline={  }
\hoffset-.3cm
\nopagenumbers
\baselineskip20pt
\centerline{  }
\vskip4truecm
\centerline{\gr  On the valuation of arithmetic--average Asian}
\vskip2pt
\centerline{\gr options: Laguerre series and Theta integrals} 
\vskip.3cm
\baselineskip14pt
\centerline{\grossrm by}
\centerline{\gross Michael Schr\"oder}
\centerline{\grossrm (Mannheim)}
\vskip.3cm
\centerline{\grossrm December 2000}
\bye